\documentclass{dcds}
\usepackage{psfig}
\usepackage{euscript}
\def\currentvolume{0}
\def\currentissue{0}
\def\currentyear{1999}
\def\month{00000}
\def\ppages{000--000}
\newtheorem{THEOREM}{Theorem}[section]
\newtheorem{LEMMA}{Lemma}[section]
\newenvironment{PROOF}{
                       \noindent{\bf Proof}.}{\qed}
\renewcommand{\theequation}{\arabic{section}.\arabic{equation}}
\newcommand{\half}{\mathrm{\frac{1}{2}}}
\newcommand{\oneovertwo}{\textstyle{\mathrm{\frac{1}{2}}}\displaystyle}
\newcommand{\threeoverfour}{\textstyle{\mathrm{\frac{3}{4}}}\displaystyle}
\newcommand{\oneoverthirtytwo}{\textstyle{\mathrm{\frac{1}{32}}}\displaystyle}
\newcommand{\eps}{\varepsilon}
\title[Vortex solutions of the CGL equation]
{Bifurcating Vortex Solutions of the \\ 
Complex Ginzburg-Landau Equation}

\author[Hans G. Kaper and Peter Tak\'{a}\v{c}]{} 

\thanks{The
work of H.G.K. is supported by
the Mathematical, Information, and
Computational Sciences Division subprogram of
the Office of Advanced Scientific Computing Research,
U.S. Department of Energy,
under Contract W-31-109-Eng-38.
The work of P.T. is supported by
the Deutsche Forschungsgemeinschaft (DFG), Germany.}


\begin{document}
\maketitle

\setcounter{page}{111}

\centerline{\scshape Hans G. Kaper}
\medskip

{\footnotesize
\centerline{Mathematics and Computer Science Division}
\centerline{Argonne National Laboratory}
\centerline{Argonne, IL 60439, USA}}

\bigskip

\centerline{\scshape Peter Tak\'{a}\v{c}}
\medskip

{\footnotesize
\centerline{Fachbereich Mathematik}
\centerline{Universit\"{a}t Rostock}
\centerline{Universit\"{a}tsplatz 1}
\centerline{D-18055 Rostock, Germany}}

\bigskip
\begin{quote}{\normalfont\fontsize{8}{10}\selectfont
{\bfseries Abstract.}
It is shown that the complex Ginzburg-Landau (CGL)
equation on the real line admits
nontrivial $2\pi$-periodic vortex solutions
that have $2n$ simple zeros (``vortices'') per period.
The vortex solutions bifurcate from
the trivial solution and inherit
their zeros from the solution
of the linearized equation.
This result rules out the possibility
that the vortices are determining nodes
for vortex solutions of the CGL equation. \par}
\end{quote}

\section{Vortex Solutions and Determining Nodes}
\setcounter{equation}{0}
In this article we investigate the bifurcation
of $2\pi$-periodic vortex solutions of the
complex Ginzburg-Landau (CGL) equation
on the real line,
\begin{equation}
   u_t
   =
   (1 + i\nu) u_{xx}
   +
   (R - (1 + i\mu) |u|^2) u ,
   \quad x \in \textbf{R}, \ t > 0 .
\label{CGL-eq}
\end{equation}
The unknown function $u$ is complex-valued;
$R$, $\mu$, and $\nu$ are given real constants.
Vortex solutions are nontrivial solutions
whose zero set consists of isolated points.
(The term ``vortex'' for a zero of $u$,
which is rather meaningless in the present context,
is borrowed from the theory of the Ginzburg-Landau equations
of superconductivity in two dimensions.
There, a zero of the complex order parameter
identifies a vortex of magnetic flux.)
The vortex solutions we are interested in
are classical solutions of the following type:
\begin{equation}
   u(x, t) = U(nx) {\rm e}^{-i\omega t} ,
   \quad x \in \textbf{R}, \ t > 0 ,
\label{u-U}
\end{equation}
where $\omega$ is a suitable real constant
that depends on $R$, $\mu$, and $\nu$,
$n$ is a fixed positive integer,
and $U$ is a $2\pi$-periodic
complex-valued $C^2$-function
that has two simple zeros per period.
(Thus $u$, which is also $2\pi$-periodic,
has $2n$ simple zeros per period.)

The investigation is motivated by the observation
that the solution of a dissipative partial
differential equation such as the CGL equation
is determined uniquely and completely
by its \textit{nodal values}---that is,
by its values at a set of determining nodes.
The concept of determining nodes was first
introduced by Foias and Temam
in the context of the Navier-Stokes equations
for viscous incompressible fluids~\cite{ft}.
These authors showed that the solution of the
two-dimensional Navier-Stokes equations
is determined uniquely and completely
by its values at a finite set of isolated points
(determining nodes).
The existence of a set of determining nodes
has since been shown for various equations,
including the CGL equation~\cite{kukavica},
the Kuramoto-Sivashinsky equation~\cite{fk},
and the Ginzburg-Landau equations
of superconductivity~\cite{kww97}.
These existence results all require that,
in some sense, the set of determining nodes be
``sufficiently dense'' in the domain,
although the cardinality of the set is unknown.
For the Navier-Stokes equation,
an upper bound of the cardinality
has been given in terms of the
physical parameters~\cite{jt93},
but it has been conjectured on the basis of the
Takens imbedding theorem~\cite{takens} that,
for dissipative partial differential equations,
the cardinality is in fact independent of
the parameters and determined entirely by
the dimensionality of the spatial domain.

By definition, if two solutions
of the CGL equation coincide at the determining nodes,
they coincide everywhere in the domain.
Since the CGL equation admits the trivial solution,
and any vortex solution coincides with
the trivial solution at the vortices,
the existence of vortex solutions
would rule out the possibility that
a solution of the CGL equation is determined
uniquely and completely by its vortices.
Indeed, an example of such a solution
satisfying the Neumann boundary conditions
on the interval $(0,1)$ was constructed by
Tak\'a\v{c}~\cite[Corollary~3.2]{takac}.
In the present work, the boundary conditions
are replaced by a condition fixing the vortices.

If $u$ is to be a vortex solution
of the type~(\ref{u-U}) with $2n$ vortices per period,
then $U$ must satisfy the nonlinear differential equation
\begin{equation}
   - U'' - U = \rho (r - |U|^2) U ,
   \quad x \in \textbf{R} ,
\label{U}
\end{equation}
where the complex constants $\rho$ and $r$
are defined in terms of $R$, $\mu$, $\nu$, and $n$,
\begin{equation}
   \rho = \frac{1 + i\mu}{(1 + i\nu) n^2} , \quad
   r = \frac{R + i\omega - (1 + i\nu) n^2}{1 + i\mu} .
\end{equation}
The problem is thus defined as a bifurcation problem,
where $\rho$ is the bifurcation parameter,
and we are interested in solutions of Eq.~(\ref{U})
that bifurcate from the trivial solution
$(r,U) = (0,0)$.

We show the following results.
First, there exist vortex solutions of the CGL equation
that have $2n$ simple zeros per period
and bifurcate from the trivial solution.
This result rules out the possibility
that the vortices are determining nodes
for vortex solutions of the CGL equation.
Second, the bifurcating vortex solutions
inherit their zeros from the solution of
the linearized equation.
The vortices that are introduced
at bifurcation are pinned
as the bifurcation parameter increases.
Moreover, numerical computations indicate
that no other zeros arise after a bifurcation.

The first result may seem to contradict
a result of Kukavica~\cite{kukavica},
who showed that the solution of the CGL equation
is completely determined by the values at two nodes,
provided these nodes are sufficiently close.
After all, by choosing $n$ sufficiently large,
we can bring the zeros of the bifurcating solution
arbitrarily close together.
However, there is no contradiction because
the upper bound on the distance
between the two determining nodes
depends on the parameters and
decreases as $n$ increases.

The linearized problem is analyzed in Section~2,
the bifurcation analysis is given in Section~3,
and numerical results are presented in Section~4.

\section{Linearized Problem}
\setcounter{equation}{0}
If Eq.~(\ref{U}) is linearized
about the trivial solution,
it reduces to
\begin{equation}
   - U'' - U = 0 ,
   \quad x \in \textbf{R} .
\label{U-lin}
\end{equation}
This equation admits $2\pi$-periodic solutions
that have two simple zeros per period.
The zeros are uniformly distributed and
separated by a distance $\pi$.

Now consider the inhomogeneous equation
\begin{equation}
   - U'' - U = f ,
   \quad x \in \textbf{R} ,
\label{U-lin-f}
\end{equation}
where $f : \textbf{R} \to \textbf{C}$
is continuous.
We claim that, under suitable conditions
on $f$, this equation admits solutions
whose zeros coincide with the zeros of
the solution of the homogeneous equation.
We make this claim precise in the following lemma
for the case where the zeros of the two solutions
coincide with the zeros of the cosine function.
Other cases are treated similarly.

\begin{LEMMA} \label{L-1}
Equation~(\ref{U-lin-f}) admits
a classical solution
that has simple zeros
at all odd multiples
of $\half\pi$
if and only if
\begin{equation}
  \int_{(k-\half)\pi}^{(k+\half)\pi} f(y) \cos y \,{\rm d}y = 0 ,
  \quad k \in \textbf{Z} .
  \label{solv-fk}
\end{equation}
If $f$ satisfies the condition~(\ref{solv-fk}),
then
\begin{equation}
  U(x) = v(x) \cos x ,
  \quad x \in \textbf{R} ,
  \label{U-v}
\end{equation}
where $v \in C^2 (\textbf{R})$ is given locally
on each interval
$[(k-\half)\pi, (k+\half)\pi]$,
$k \in \textbf{Z}$,
by the expression
\begin{equation}
  v(x)
  =
  v((k-\oneovertwo)\pi)
  +
  \int_{(k-\half)\pi}^{(k+\half)\pi}
  f(y)
  \frac{g(x,y)}{\cos x}
  \,{\rm d}y .
  \label{vk}
\end{equation}
The kernel $g$ is independent of $k$,
\begin{equation}
  g(x,y) = \left\{
  \begin{array}{ll}
  \cos x \sin y & \mbox{if }\, y \leq x , \\
  \sin x \cos y & \mbox{if }\, y \geq x .
  \end{array}
  \right.
  \label{g}
\end{equation}
\end{LEMMA}

\begin{PROOF}
Let $f: \textbf{R} \to \textbf{C}$
be a given continuous function.
If we look for a solution $U$
of Eq.~(\ref{U-lin-f})
of the form~(\ref{U-v}),
then $v$ must satisfy
the degenerate differential equation
\begin{equation}
  - v'' \cos x + 2 v' \sin x = f
  \label{vf-eq}
\end{equation}
for all $x \not= (k+\half) \pi$;
moreover, $v$ must remain bounded
near the points $(k+\half) \pi$,
for all $k \in \textbf{Z}$.

Equation~(\ref{vf-eq}) can be integrated
locally on any interval
$((k-\half) \pi, (k+\half) \pi)$.
In fact, after multiplying both sides
of the equation by $\cos x$, we have
\begin{equation}
  - (v' \cos^2 x)' = f(x) \cos x .
  \label{v'-eq}
\end{equation}
If $v_k$ is the local representation of $v$
on $((k-\half) \pi , (k+\half) \pi )$,
then the integration yields
\[
  v_k'(x) \cos^2 x
  =
  v_k' (k\pi)
  - \int_{k\pi}^x f(y) \cos y \, {\rm d}y ,
  \quad (k-\oneovertwo) \pi < x < (k+\oneovertwo) \pi .
\]
For $v_k'$ to remain bounded near the endpoints
$(k\pm\half) \pi$,
it is necessary and sufficient that
\[
  v_k' (k\pi)
  =
  \int_{k\pi}^{(k \pm \half) \pi}
  f(y) \cos y \, {\rm d}y ,
\]
so $f$ must satisfy the solvability
condition~(\ref{solv-fk}).

If $f$ satisfies the condition~(\ref{solv-fk}),
then
\[
  v_k'(x)
  =
  - \frac{1}{\cos^2 x}
  \int_{(k\pm\half)\pi}^x f(y) \cos y \, {\rm d}y ,
  \quad
  (k-\oneovertwo) \pi < x < (k+\oneovertwo) \pi ,
\]
and
$v_k'((k\pm\half)\pi) = \mp \half f((k\pm\half)\pi)$.
The expression~(\ref{vk}) follows upon integration.
\end{PROOF}

While Eq.~(\ref{vk}) gives a local representation
of $v$ on each interval
$[(k-\half)\pi$, $(k+\half)\pi]$,
there also exists a global representation
that is valid on the entire real line.
First, observe that
\begin{equation}
  v((k+\oneovertwo)\pi)
  =
  v((k-\oneovertwo)\pi)
  +
  \int_{(k-\half)\pi}^{(k+\half)\pi} f(y) \sin y \,{\rm d}y ,
  \quad k \in \textbf{Z} .
  \label{v-cont}
\end{equation}
Repeated application of this recurrence relation
yields an expression for $v((k-\half)\pi)$
in terms of $v(-\half\pi)$,
\[
  v ((k-\oneovertwo)\pi)
  =
  v (-\oneovertwo\pi)
  +
  \int_{-\half\pi}^{(k-\half)\pi} f(y) \sin y \, {\rm d}y ,
  \quad k \in \textbf{Z} .
\]
Furthermore, because $f$ satisfies~(\ref{solv-fk}),
\[
  \int_x^{(k+\half)\pi} f(y) \cos y \,{\rm d}y
  =
  \int_x^{\half\pi} f(y) \cos y \,{\rm d}y ,
  \quad k \in \textbf{Z} .
\]
Thus, $v$ is represented globally by the expression
\begin{equation}
  v(x)
  =
  v (-\oneovertwo\pi)
  + \int_{-\half\pi}^x
  f(y) \sin y \, {\rm d}y
  + \frac{\sin x}{\cos x} \int_x^{\half\pi}
  f(y) \cos y \, {\rm d}y ,
  \quad x \in \textbf{R} .
  \label{v}
\end{equation}

\section{Bifurcation Analysis}
\setcounter{equation}{0}
We now proceed to the bifurcation analysis.
We recall that we wish to find solutions
of Eq.~(\ref{U}) that are $2\pi$-periodic
and have two simple zeros per period.
In fact, we will try to find solutions
whose zeros coincide with the zeros of $\cos x$---the
solution of the linearized equation.

We use the results of the preceding section,
substituting for $f$ the expression
in the right member of Eq.~(\ref{U}).
Taking $U$ to be of the form
(cf.~\cite[Eq.~(3.19)]{takac})
\begin{equation}
  U(x) = v(x) \cos x ,
  \quad x \in \textbf{R} ,
  \label{Uv}
\end{equation}
we replace the original problem
by a bifurcation problem for $(r,v)$
in a neighborhood of
$(r,v) = (0,0) \in \textbf{C} \times C^2(\textbf{R})$.

We infer from Lemma~\ref{L-1}
that the bifurcation analysis
can  performed locally
on any of the intervals
$[(k-\half)\pi, (k+\half)\pi]$,
$k \in \textbf{Z}$.
Hence, it suffices
to consider the function $v$ on
the interval $[-\half\pi,\half\pi]$,
which we denote by $J$ from now on.
According to Eq.~(\ref{vk}),
$v$ must satisfy the following
integral equation on $J$:
\begin{equation}
  v(x)
  =
  v(-\oneovertwo\pi)
  +
  \int_J f(y) \frac{g(x,y)}{\cos x} \, {\rm d}y ,
  \quad x \in J ,
  \label{v-f}
\end{equation}
where $g$ is defined in Eq.~(\ref{g})
and $f$ is given in terms of $v$,
\begin{equation}
  f(x)
  = \rho (r - |v|^2 \cos^2 x) v \cos x ,
  \quad v \equiv v(x) ,
  \quad x \in J .
  \label{f-v}
\end{equation}
The function $f$ must satisfy
the condition~(\ref{solv-fk})
for $k=0$.
With $f$ given by Eq.~(\ref{f-v}),
the latter translates into a
relation between $r$ and $v$,
\begin{equation}
  r
  \int_J
  v(y) \cos^2 y \,{\rm d}y
  =
  \int_J
  | v(y) |^2 v(y) \cos^4 y \,{\rm d}y .
  \label{r-v}
\end{equation}
If we take this as the definition of $r$,
then we have reduced the bifurcation problem
to a problem for $v$ in the neighborhood
of $v = 0 \in C(J)$.

We employ the Lyapunov-Schmidt reduction method
in much the same way as in
\cite[Proof of Theorem~3.1]{takac}.
Let the projection $P : C(J) \to C(J)$ be defined by
\begin{equation}
   P u (x)
   =
  \frac{2}{\pi}
  \int_J u(y) \cos^2 y \,{\rm d}y ,
  \quad u \in C(J) ,
  \quad x \in J ,
  \label{P}
\end{equation}
and its complement $P': C(J) \to C(J)$ by
$P'= I-P$.
($I$ is the identity operator in $C(J)$.)
The pair $(P, P')$ decomposes the space $C(J)$.
Note that $Pu$ is a complex constant-valued function,
so we may identify $PC(J)$ with the complex plane.
Note also that $P1=1$.

Let $C_0(J)$ denote the closed subspace of $C(J)$
consisting of all elements $f \in C(J)$
that satisfy the condition~(\ref{solv-fk})
for $k=0$.
For any $f \in C_0 (I)$, we define
$v \in C(I)$ by the relation~(\ref{v-f});
its projection $Pv$ is
\begin{equation}
  Pv (x)
  =
  v(-\oneovertwo\pi)
  +
  \int_J
  \left(
  \frac{2}{\pi}
  \int_J g(z,y) \cos z \,{\rm d}z
  \right)
  f(y) \,{\rm d}y ,
  \quad x \in J .
  \label{Pv}
\end{equation}
We set $Pv = \eps$ and scale $P'v$ by $\eps$,
putting $P'v = \eps w$.
Thus,
\begin{equation}
  v = \eps (1+w) ,
  \quad \eps \in \textbf{C} ,
  \quad w \in P'C(J) .
  \label{w}
\end{equation}
The mapping $f \mapsto \eps w$
defines a linear operator $L$
from $C_0(J)$ into $P'C(J)$,
\begin{equation}
  L f
  =
  \eps w ,
  \quad f \in C_0(J) .
  \label{w-f}
\end{equation}
Since $\eps w = v - Pv$,
the expression for $Lf$ is readily found
from Eqs.~(\ref{v-f}) and~(\ref{Pv}),
\begin{equation}
  (Lf)(x) = \int_J
  \left(
  \frac{g(x,y)}{\cos x}
  - \frac{2}{\pi} \int_J g(z,y) \cos z \,{\rm d}z
  \right)
  f(y) \,{\rm d}y ,
  \quad x \in J .
  \label{Lf}
\end{equation}

\begin{LEMMA} \label{L-f}
The linear operator $L : C_0(J) \to C(J)$
defined in Eq.~(\ref{w-f}) is bounded,
\begin{equation}
   \| Lf \|_\infty \leq 3\pi \| f \|_\infty ,
   \quad f \in C_0(I) .
\end{equation}
\end{LEMMA}

\begin{PROOF}
Since $|g(x,y)| \leq 1$,
it is certainly true that
\begin{equation}
  \left|
  \int_J
  \frac{2}{\pi}
  \int_J g(z,y) \cos z
  \,{\rm d}z\
  f(y) \,{\rm d}y
  \right|
  \leq
  \frac{2}{\pi} |J|^2
  \|f\|_\infty
  =
  2\pi \|f\|_\infty ,
  \quad x \in J .
  \label{est1}
\end{equation}
To estimate the remaining integral
in Eq.~(\ref{Lf}),
we distinguish between $x \geq 0$
and $x \leq 0$.

Suppose $x \geq 0$.
Then
\[
  \int_J
  \frac{g(x,y)}{\cos x}
  f(y) \, {\rm d}y
  =
  \int_{-\half\pi}^x f(y) \sin y \,{\rm d}y
  + \sin x 
  \int_x^{\half\pi} f(y) \frac{\cos y}{\cos x} \,{\rm d}y .
\]
The first term is estimated trivially;
its modulus is less than or equal to
$(x+\half\pi) \|f\|_\infty$.
In the second term, we use
the fact that
$0 \leq \cos y / \cos x \leq 1$
for all $0 \leq x < y \leq \half\pi$;
the modulus of this term is less than
$(\half\pi - x) \| f \|_\infty$.
Together, these two inequalities
give the estimate
\begin{equation}
  \left|
  \int_J
  \frac{g(x,y)}{\cos x}
  f(y) \, {\rm d}y
  \right|
  \leq
  \pi \| f \|_\infty ,
  \quad x \in J ,
  \quad x \geq 0 .
  \label{est2-pos}
\end{equation}
Now suppose $x \leq 0$.
Then we start from the expression
\[
  \int_J
  \frac{g(x,y)}{\cos x}
  f(y) \, {\rm d}y
  =
  \int_{-\half\pi}^x f(y) \sin y \,{\rm d}y
  + \sin x 
  \int_{-\half\pi}^x f(y) \frac{\cos y}{\cos x} \,{\rm d}y
\]
and find, similarly,
\begin{equation}
  \left|
  \int_J
  \frac{g(x,y)}{\cos x}
  f(y) \, {\rm d}y
  \right|
  \leq
  \pi \| f \|_\infty ,
  \quad x \in J ,
  \quad x \leq 0 .
  \label{est2-neg}
\end{equation}
Together, the inequalities~(\ref{est2-pos})
and~(\ref{est2-neg}) give the estimate
\begin{equation}
  \left|
  \int_J
  \frac{g(x,y)}{\cos x}
  f(y) \, {\rm d}y
  \right|
  \leq
  \pi \| f \|_\infty ,
  \quad x \in J .
  \label{est2}
\end{equation}
The statement of the lemma follows from
Eqs.~(\ref{Lf}), (\ref{est1}), and~(\ref{est2}).
\end{PROOF}

The integral in the left member of Eq.~(\ref{r-v})
is equal to
$\half\pi Pv$, where $Pv = \eps$,
so the condition~(\ref{r-v}),
which we use to define $r$ in terms of $v$,
reduces to
\begin{equation}
  r =
  \frac{2}{\eps\pi}
  \int_J |v(y)|^2 v(y) \cos^4 y \,{\rm d}y .
\label{r-w}
\end{equation}
When we insert this expression into Eq.~(\ref{f-v})
and make the substitution $v = \eps(1+w)$,
we obtain a relation between $f$ and $w$,
\begin{equation}
   f
   =
   \eps |\eps|^2 F(w) , \quad w \in P'C(J) ,
   \label{f-F}
\end{equation}
where $F: P'C(J) \to C_0(J)$ is the following nonlinear map:
\[
   [F(w)](x)
   =
   \rho \left(
   \frac{2}{\pi}
   \int_J |1+w(y)|^2 (1+w(y)) \cos^4 y \,{\rm d}y
   - |1+w(x)|^2 \cos^2 x
   \right)
\]
\begin{equation}
  \times
  (1+w(x)) \cos x ,
  \quad x\in J ,
  \, w\in P'C(J) .
\label{G-f}
\end{equation}
Combining Eqs.~(\ref{w-f}) and~(\ref{f-F}),
we obtain an equation for $w$ in $P'C(J)$,
\begin{equation}
   w = T_\eps(w)
   = | \eps |^2 L ( F (w) ) .
\label{w-eq}
\end{equation}
We wish to solve this equation
using the Banach contraction
principle~\cite[Theorem~7.1]{deimling}.
We already know that $L$ is bounded
from $C_0(J)$ into $P'C(J)$;
the following lemma gives
the necessary estimates for $F$.

Let ${\EuScript B}_\sigma$ denote
the closed ball of radius $\sigma$
($\sigma > 0$)
centered at the origin in $P'C(J)$,
\begin{equation}
   {\EuScript B}_{\sigma}
   = \{ w \in P'C(J) : \|w\|_\infty \leq \sigma \} .
\end{equation}

\begin{LEMMA} \label{L-F}
The nonlinear map
$F : P'C(J) \to C_0(J)$
defined in Eq.~(\ref{f-F})
is bounded and Lipschitz continuous,
\begin{equation}
   \| F(w) \|_\infty
   \leq
   |\rho| (2+\sigma) (1+\sigma)^3 ,
  \quad w \in {\EuScript B}_\sigma ,
\label{Fw}
\end{equation}
\begin{equation}
   \| F(w_1) - F(w_2)\|_\infty
   \leq
   3 |\rho| (2+\sigma) (1+\sigma)^2
   \| w_1 - w_2 \|_\infty ,
  \quad w_1, w_2 \in {\EuScript B}_\sigma .
\label{F1-F2}
\end{equation}
\end{LEMMA}

\begin{PROOF}
If $w \in {\EuScript B}_\sigma$, then
\[
   | [F(w)](x) |
   \leq
   |\rho|
   \left(
   \frac{2}{\pi} 
   (1 + \sigma)^3
   \int_J \cos^4 y \,{\rm d}y
   +
   (1 + \sigma)^2
   \right)
   (1 + \sigma) ,
   \quad x \in J .
\] 
Because $(2/\pi) \int_J \cos^4 y \, {\rm d}y = \frac{3}{4} < 1$,
the estimate~(\ref{Fw}) follows.

If $w_1, w_2 \in {\EuScript B}_\sigma$, then
\begin{eqnarray*}
&  \left| [F(w_1)] (x) - [F(w_2)] (x) \right|
   \leq
   |\rho|
   \left| s_1 (1+w_1 (x)) - s_2 (1+w_2 (x)) \right|
\\
&  +
   |\rho|
   \left|
   |1+w_1 (x)|^2 (1+w_1 (x)) - |1+w_2 (x)|^2 (1+w_2 (x))
   \right| ,
\end{eqnarray*}
where we have used the abbreviations
\[
   s_j = \frac{2}{\pi}
   \int_J |1+w_j(y)|^2 (1+w_j(y)) \cos^4 y \,{\rm d}y , \quad
   j=1,2 .
\]
Adding and subtracting terms, we see that
\begin{eqnarray*}
&  \left|
   |1+w_1|^2 (1+w_1) - |1+w_2|^2 (1+w_2)
   \right|
\\
&  =
   \left|
   ( |1+w_1|^2 + |1+w_2|^2 ) (w_1 - w_2)
   + (1 + w_1) (1 + w_2) (\overline{w}_1 - \overline{w}_2)
   \right|
\\
&  \leq
   3 (1+\sigma)^2 \| w_1 - w_2 \|_\infty .
\end{eqnarray*}
Furthermore,
\begin{eqnarray*}
&  \left|
   s_1 (1+w_1) - s_2 (1+w_2)
   \right|
   =
   \left|
   (1+w_1) (s_1 - s_2) + s_2 (w_1 - w_2)
   \right|
\\
&  \leq
   (1+\sigma) | s_1 - s_2 | + |s_2| | w_1 - w_2 | .
\end{eqnarray*}
One readily verifies that
\[
   |s_1 - s_2| \leq
   \frac{6}{\pi}
   (1+\sigma)^2
   \left( \int_J \cos^4 y \,{\rm d}y \right)
   \| w_1 - w_2 \|_\infty
   =
   \frac{9}{4} (1+\sigma)^2 \| w_1 - w_2 \|_\infty
\]
and
\[
   |s_2|
   \leq
   \frac{2}{\pi}
   (1+\sigma)^3
   \left( \int_J \cos^4 y \,{\rm d}y \right)
   =
   \frac{3}{4} (1+\sigma)^3 ,
\]
so
\[
   \left|
   s_1 (1+w_1) - s_2 (1+w_2)
   \right|
   \leq
   3 (1+\sigma)^3 \| w_1 - w_2 \|_\infty .
\]
The inequality~(\ref{F1-F2}) follows.
\end{PROOF}

We are ready to prove the desired bifurcation result.
Let the set $\Gamma$ be defined by
\[
  \Gamma
  =
  \{
  (r,U) \in \textbf{C} \times C^2 (\textbf{R}) :
  (r,U) \mbox{ satisfies Eq.~(\ref{U})} ;
\]
\begin{equation}
  U(x) = v(x) \cos x , \, x \in \textbf{R} ;
  \ v \in C(\textbf{R}) 
  \ v \mbox{ bounded}
  \}
\label{Gamma}
\end{equation}

\begin{THEOREM} \label{T1}
The point
$(0,0) \in \textbf{C} \times C^2 (\textbf{R})$
is a bifurcation point for Eq.~(\ref{U}).
There exists an open neighborhood ${\EuScript O}$
of $(0,0)$ in $\textbf{C} \times C^2 (\textbf{R})$
and a positive constant $\delta$ such that
the set $\Gamma \cap {\EuScript O}$
coincides with the set of all
$(r,U) \in \textbf{C} \times C^2 (\textbf{R})$
having the following representation:
\begin{equation}
   r
   =
   \threeoverfour |\eps|^2
   \left(1 + |\eps|^2 \varphi (|\eps|^2 \right) ,
\label{r-rep}
\end{equation}
\begin{equation}
   U(x)
   = 
   \eps
   (1 + |\eps|^2 \Phi (|\eps|^2, x) )
   \cos x ,
   \quad x \in\textbf{R} , 
\label{U-rep}
\end{equation}
where $\eps$ is an arbitrary complex parameter
with $0 < | \eps |^2 < \delta$,
and
$\varphi : [0,\delta) \to \textbf{C}$
and
$\Phi : [0,\delta) \times \textbf{R} \to \textbf{C}$
are continuous functions satisfying the following conditions:

\noindent{\rm (i)}
$(r,U) \in \Gamma$,

\noindent{\rm (ii)}
$\Phi (s, \cdot) \in C^2 (\textbf{R})$
for every $s \in (0, \delta)$
and $\int_{\textbf{R}} \Phi (s,x) \cos^2 x \,{\rm d}x = 0$, and

\noindent{\rm (iii)}
the real and imaginary parts of $\varphi$ and $\Phi$
are real-analytic functions of their arguments.
\end{THEOREM}

\begin{PROOF}
Following the steps outlined in the preceding analysis,
we reduce the bifurcation problem to a problem for $w$
in the neighborhood of $w = 0 \in P'C(J)$.
This function $w$ must be a fixed point of
the operator $T_\eps$ defined in Eq.~(\ref{w-eq}).
Once $w$ has been found, we define $v$ in terms
of $w$ by means of Eq.~(\ref{w}) and
$(r,U)$ in terms of $v$ by means of
Eqs.~(\ref{r-w}) and~(\ref{Uv}).

From Lemmas~\ref{L-f} and~\ref{L-F}
we obtain
\[
   \| L (F (w)) \|_\infty
   \leq
   3\pi \| F (w) \|_\infty
   \leq
   3\pi |\rho| (2+\sigma) (1+\sigma)^3 ,
  \quad w \in {\EuScript B}_\sigma ,
\]
so $T_\eps = |\eps|^2 L F$
maps ${\EuScript B}_\sigma$ into itself whenever
\[
   | \eps |^2
   <
   \frac{\sigma}{3\pi|\rho| (2 + \sigma)(1 + \sigma)^3 } .
\]
Furthermore,
\[
   \| T_\eps (w_1 - w_2) \|_\infty
   \leq
   9\pi
   | \eps |^2
   |\rho| (2+\sigma)(1+\sigma)^2
   \| w_1 - w_2 \|_\infty ,
  \quad w_1, w_2 \in {\EuScript B}_\sigma ,
\]
so $T_\eps $ is a contraction if
\begin{equation}
   | \eps |^2
   <
   \frac{1}{9\pi|\rho| (2+\sigma)(1+\sigma)^2} .
\end{equation}
Hence, if we define
\begin{equation}
   \delta \equiv
   \delta (\sigma) =
   \frac{1}{3\pi|\rho| (2 + \sigma) (1+\sigma)^2}
   \min
   \left\{
   \frac{\sigma}{1 + \sigma} ,
   \frac{1}{3}
   \right\} ,
\end{equation}
then $T_\eps $ is a contractive mapping of
${\EuScript B}_\sigma$ into itself
for every $\varepsilon \in \textbf{C}$
satisfying $0 < |\varepsilon|^2 < \delta$.
Consequently, $T_\varepsilon$ has
a unique fixed point in ${\EuScript B}_\sigma$,
which can be found by iteration.
The lowest-order approximation $w=0$,
which corresponds to $v=\eps$,
gives $r=\frac{3}{4} |\varepsilon|^2$
and $U(x) = \eps \cos x$.

The statements of the theorem
follow from the implicit function
theorems~\cite[Theorems~15.1 and~15.3]{deimling}.
\end{PROOF}

Theorem~\ref{T1} implies that the CGL equation
admits $2\pi$-periodic vortex solutions $u$,
which bifurcate from the trivial solution;
these vortex solutions have $2n$ zeros
(``vortices'') per period; and
the vortices are located at the zeros
of the cosine function, which is the solution
of the linearized equation in the neighborhood
of the bifurcation point.

The conditions~(i)--(iii),
together with the representations~(\ref{r-rep}) and~(\ref{U-rep}),
determine $\eps$, $\varphi$, and $\Phi$ uniquely.

The representations~(\ref{r-rep}) and~(\ref{U-rep}) show that
we have a supercritical pitchfork bifurcation from $(0,0)$.
Further terms in the representations~(\ref{r-rep}) and~(\ref{U-rep})
can be computed in a standard manner,
\begin{equation}
   r
   =
   \threeoverfour
   |\eps|^2
   \left(
   1
   -
   \oneoverthirtytwo\rho|\eps|^2
   +
   O( |\eps|^4 )
   \right) ,
\label{r-as}
\end{equation}
\begin{equation}
   U(x)
   = 
   \eps
   \cos x
   \left(
   1 -
   \frac{\rho|\eps|^2}{32}
   \frac{\cos 3x}{\cos x}
   +
   \left(\frac{\rho|\eps|^2}{32}\right)^2
   \frac{3\cos 3x + \cos 5x}{\cos x}
   +
   O( |\eps|^6 )
   \right) .
\label{U-as}
\end{equation}
Furthermore,
\begin{equation}
   R = n^2 + \threeoverfour
   |\eps|^2
   \left(
   1
   -
   \frac{1 - \mu^2 + 2\mu\nu}{32 n^2 (1 + \nu^2)}
   |\eps|^2
   +
   O( |\eps|^4 )
   \right) ,
\label{R-as}
\end{equation}
and
\begin{equation}
   \omega = \nu n^2 + \threeoverfour
   |\eps|^2
   \left(
   \mu
   -
   \frac{\mu^2\nu + 2\mu - \nu}{32 n^2 (1 + \nu^2)}
   |\eps|^2
   +
   O( |\eps|^4 )
   \right) .
\label{omega-as}
\end{equation}
In particular,
$\omega = \mu R + (\nu - \mu) n^2 + O((R - n^2)^2)$.

\section{Numerical Results}
\setcounter{equation}{0}
The results of the preceding bifurcation analysis
are supported by the results of numerical computations.
(These computations were performed by Michael Levine,
participant in the 1998 Energy Research
Undergraduate Laboratory Fellowship program
at Argonne National Laboratory.)

Three numerical methods were applied.
The first method was a fixed-point iteration
based on Eq.~(\ref{w-eq}).
(Observe that the only parameter in Eq.~(\ref{w-eq})
is $\rho |\eps|^2$;
without loss of generality, we may take $\eps = 1$.)
The method converged for $\rho$ in the rectangle
$[-3.5, 3.5] \times [0, 1.5]$.
The bifurcating solutions were found to be
very close to the solutions of the
linearized equation.
Next, a shooting method was applied to Eq.~(\ref{U}).
The method yielded bifurcating solutions for $\rho$
in discs centered at the origin with
radii up to 9.
In a third method, a finite-difference
method was applied to Eq.~(\ref{U}), and
the resulting system of linear equations
was solved directly.
This method gave results for $\rho$
in discs centered at the origin
with radii up to 200.

None of the bifurcating solutions
had any additional zeros.
The bifurcating solutions were all symmetric
with respect to the origin.
For values of $\rho$ close to the imaginary axis,
additional asymmetric solutions were found
that bifurcated from the symmetric ones.
These bifurcations occurred multiple times
as $|\rho|$ was increased along rays emanating
from the origin, and we conjecture that they
occur infinitely often.

The properties of the bifurcating solutions
are summarized in Figs.~1 and~2.

\vspace{.5in}
\begin{figure}[htb]
\begin{center}
\vspace{-5ex}
\mbox{\psfig{file=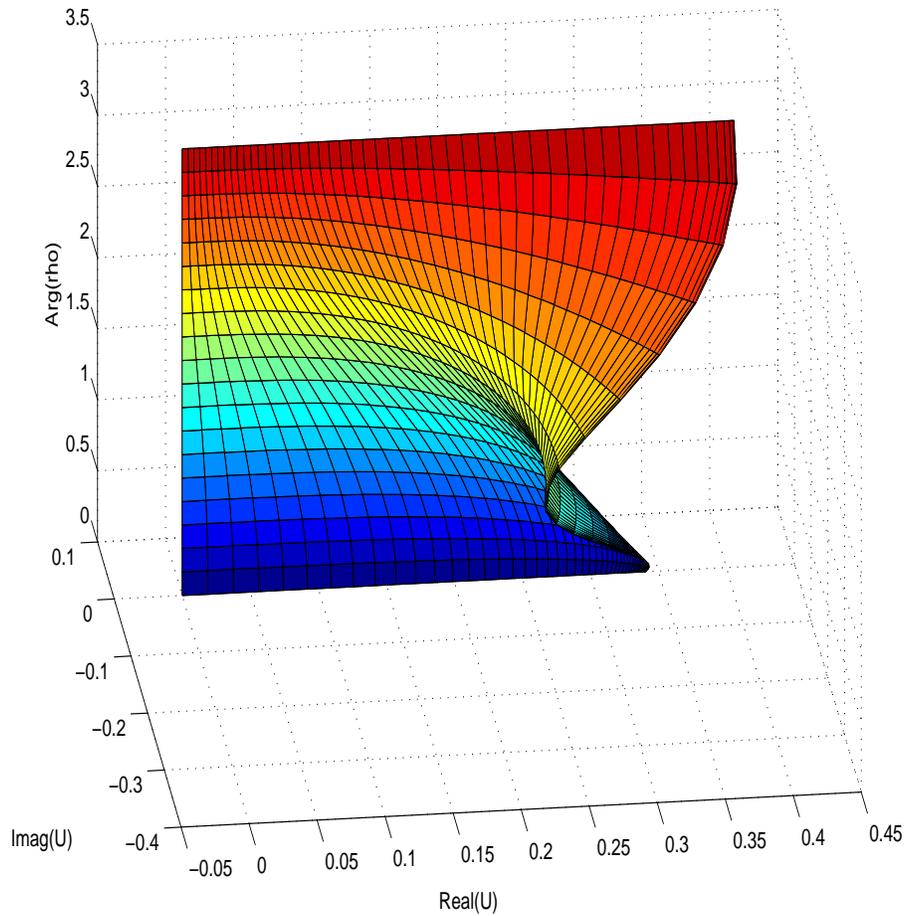,height=5in,width=4.7in}}
\vspace{-2ex}
\caption{
Bifurcating solutions $U$ of the CGL equation
as a function of $\arg(\rho)$ for a fixed value of $|\rho|$.}
\end{center}
\end{figure}

\newpage
\vspace*{.2in}
\begin{figure}[htb]
\begin{center}
\vspace{-5ex}
\mbox{\psfig{file=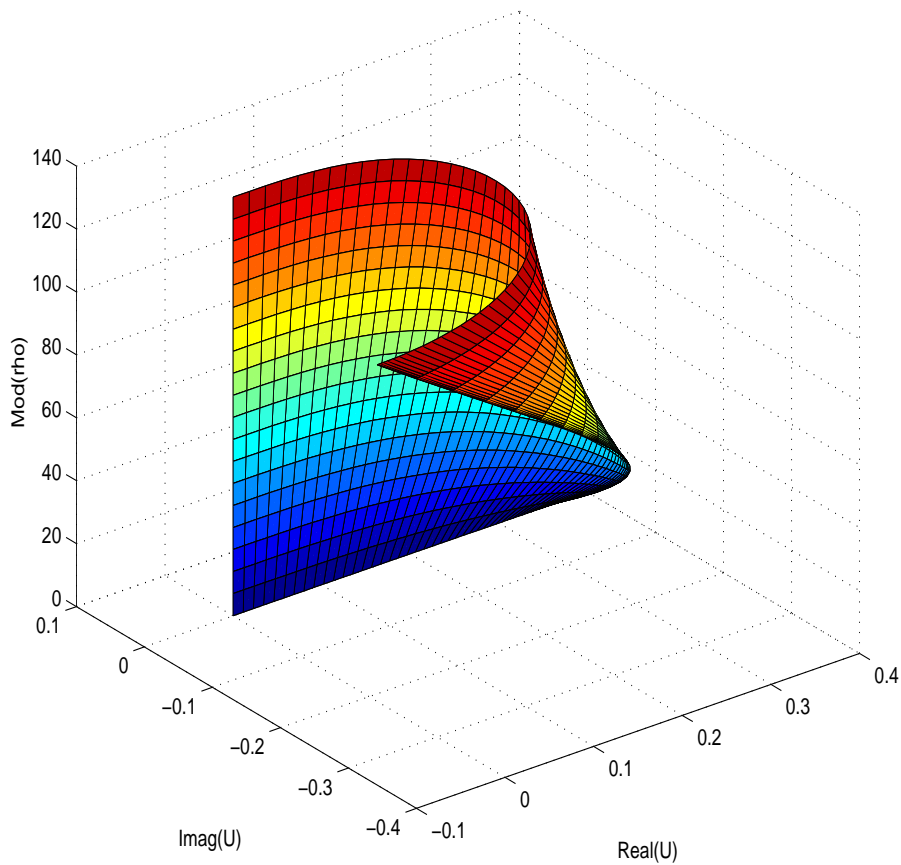,height=4.72in,width=4.7in}}
\vspace{-2ex}
\caption{
Bifurcating solutions $U$ of the CGL equation
as a function of $|\rho|$ for a fixed value of $\arg(\rho)$.}
\end{center}
\end{figure}

\noindent
\textit{E-mail addresses}:
\texttt {kaper@mcs.anl.gov, peter.takac@mathematik.uni-rostock.de}

\begin{thebibliography}{99}

\bibitem{deimling}
K.~Deimling,
\textit{Nonlinear Functional Analysis},
Springer-Verlag, Berlin, 1985.

\bibitem{fk}
C.~Foias and I.~Kukavica,
\textsl{Determining nodes for the Kuramoto-Sivashinsky equation},
J.~Dynam.\ Diff.\ Eq.\ {\bf 7} (1995), 365--373.

\bibitem{ft}
C.~Foias and R.~Temam,
\textsl{Determination of the solutions of the
Navier-Stokes equations by a set of nodal values},
Math.\ Comp. {\bf 43} (1984), 117--133.

\bibitem{jt93}
D.~A.~Jones and E.~S.~Titi,
\textsl{Upper bounds on the number of determining
modes, nodes, and volume elements for the Navier-Stokes equations},
Indiana Univ.\ Math.\ J.\ {\bf 42} (1993), 875--887.

\bibitem{kww97}
H.~G.~Kaper, B.~Wang, and S.~Wang,
\textsl{Determining nodes for the Ginzburg-Landau equations
of superconductivity},
Discrete and Continuous Dynamical Systems
{\bf 4} (1998), 205--224.

\bibitem{kukavica}
I.~Kukavica,
\textsl{On the number of determining nodes for the
Ginzburg-Landau equation},
Nonlinearity {\bf 5} (1992), 997--1006.

\bibitem{takac}
P.~Tak\'a\v{c},
\textsl{Invariant 2-tori in the time\--dependent
Ginzburg\--Landau equation},
Nonlinearity {\bf 5} (1992), 289--321.

\bibitem{takens}
F.~Takens,
\textsl{Detecting strange attractors in turbulence}.
In: D.~A.~Raud and L.-S.~Young (eds.),
Lecture Notes in Math., Vol. 898, Springer-Verlag,
New York, pp.~366--381.

\end{thebibliography}
\end{document}